\input amstex
\documentstyle {amsppt}
\pageheight{50.5pc}
\pagewidth{32pc}
\topmatter
\title{On stochastic continuity of generalized diffusion
processes constructed as the strong solution to an SDE}
\endtitle
\author
 Ludmila L. Zaitseva
\endauthor
\address
Kyiv Taras Shevchenko university, Volodymyrska 64, Kyiv, Ukraine,
01033
\endaddress
\abstract The comparison theorem for skew Brownian motions is
proved. As the corollary we get the estimate on ${\Cal
L}_1-$distance between two skew Brownian motions started from
different points. Using this result we prove the continuous
dependence on starting point of one class of generalized diffusion
processes constructed as the strong solution to an SDE.
\endabstract
\email
zaitseva\@mail.univ.kiev.ua
\endemail
\endtopmatter
\document
\rightheadtext{On stochastic continuity of generalized diffusion
processes}
\leftheadtext{Zaitseva L.L.}

\head Introduction\endhead

The problem we consider in this paper is estimation of the distance
between two strong solutions to SDE with singular coefficients.
Considered processes belong to the class of generalized diffusion
processes, their drift vectors and diffusion matrices include
delta-function concentrated on a hyperplane.

The class of generalized diffusion processes was introduced by
Portenko M.I. (see \cite{1}). One of the most known representative
of this class is skew Brownian motion. Firstly it appears in
monograph by It\^o K. and McKean H.P. (see \cite{2}, Section 4.2,
Problem 1), then it is constructed by Portenko M.I. as a
generalized diffusion processes in \cite{1} and by Walsh J.B.
\cite{3} in the terms of its scale function and speed measure.
Harrison J.M. and Shepp L.A. prove (see \cite{4}) that skew
Brownian motion can be constructed as the strong solution to an
SDE.

This fact allows us to consider a family of such processes indexed
by starting points or skewing parameters on the same probability
space. It occurs that this family has new properties in comparison
with solutions of standard SDE's or processes with reflection. For
example, using the It\^o formula for these classic processes one
can obtain the estimate for ${\Cal L}_p-$distance between such
processes starting from different points for all $p\geq 1.$
However, one can obtain as the corollary of the results of Burdzy
K. and Kaspi H. (see \cite{5}) that skew Brownian motion is not
continuous function of the starting point. This means that there
does not exist good estimate on ${\Cal L}_p-$distance between two
skew Brownian motions starting from different points for $p>2.$
Therefore the estimation of distance between two skew Brownian
motions and, moreover, between two strong solutions to SDE with
singular coefficients is non-trivial problem which demands new
technique to deal with. We will use the results of this paper in
our next paper devoted to the Markov property of solutions to SDE
with singular coefficients.

The plan of the paper is the following one. In Section 1 we prove
the comparison theorem for skew Brownian motions, the simple
corollary of this theorem will be the estimate for ${\Cal
L}_1-$distance between two skew Brownian motions. We use
approximation approach for proving this result. It is known that
skew Brownian motion can be constructed (see \cite{6}, p.111) as
the weak limit of an appropriate sequence of diffusion processes.
We prove that a pair of skew Brownian motions constructed as the
functional of the same Wiener process can be approximated by a
pair of diffusion process. This result together with the known
comparison theorem (see \cite {7}, Section VI, Theorem 1.1) for
diffusion processes gives us the required result. In Section 2 we
use the estimate for distance between skew Brownian motions to
prove that solutions to SDE with singular coefficients depend
continuously on starting point.

\head 1. The comparison theorem for skew Brownian motions\endhead

Let $(\Omega,{\Cal F},{\Bbb P})$ be a probability space. Consider
one-dimensional Wiener process $\{w(t)\}$ started from $0$ and
filtration ${\Cal F}_t^w=\sigma  \left\{w(u),0\leq u\leq t
\right\},t\geq 0.$ For given parameters $q\in[-1,1]$ and $x_0\in\Re$
one can construct (see \cite{4}) a pair of $\{{\Cal F}_t^w\}$-adapted
processes $\{(x(t),\eta_t)\}$ such that $\{\eta_t\}$ is the local
time in $0$ for $\{x(t)\}$ and the equality $x(t)=x_0+q\eta_t+w(t),
t\geq 0$ is true. The process $\{x(t)\}$ is called skew Brownian
motion.

For $i=1,2,$ for given parameters $q_i\in [-1,1]$ and $x^i_0\in\Re$
let us construct a pair of skew Brownian motions as the functional of
the one Wiener process $\{w(t)\}$ in such a way:
$x^i(t)=x^i_0+q_i\eta_t^i+w(t), t\geq 0.$

\proclaim{Theorem 1} (Comparison theorem for skew Brownian motions.)
Let $q_1,q_2\in(-1,1)$ and
$$
1)\ x_0^1\leq x^2_0;\quad 2)\ q_1\leq q_2.
$$
Then $x_1(t)\leq x_2(t)$ for all $t\geq 0$ a.s. \endproclaim

The proof is based on an appropriate approximation procedure for
the processes $\{x_1(t)\},$ $\{x_2(t)\}$ by diffusion processes.
In a sequel we denote by the symbol $\mathop{\longrightarrow}
\limits^W$ the weak convergence of sequences of distributions of the
processes, considered as random elements of $C([0,+\infty),{\Bbb
X}),$ where ${\Bbb X}$ is equal $\Re^1,\Re^2$ or $\Re^3$ according to
context. The following limit theorem for one skew Brownian motion is
known (see \cite {6}, p.111).

\proclaim{Proposition 1} Consider a sequence of diffusion
processes in $\Re:$
$$x_n(t)=x_0+\int_0^ta_n(x_n(\tau))d\tau+w(t),\ t\geq 0, n\geq
1,$$ where $a_n(x)=na(nx), x\in\Re, n\geq 1,$ the function
$a:\Re\to\Re$ satisfies conditions $1)\int_\Re|a(x)|dx<\infty,
2)|a(x)-a(y)|<K|x-y|, x,y\in\Re,$ for some $K>0.$ Then $x_n(\cdot)
\mathop{\longrightarrow}\limits^W x(\cdot),n\to+\infty,$ where the
process $\{x(t)\}$ is skew Brownian motion with skewing parameter
$q=\,\hbox{{\rm th}}\, A,A=\int_\Re a(x)dx.$\endproclaim

The idea of proof of Theorem 1 is to approximate a pair of skew
Brownian motions by a pair of diffusion processes and then apply
the comparison theorem for diffusion processes. We arrange the
approximation procedure in two steps.

\proclaim {Lemma 1} In a situation of Proposition 1 we have
$$
\vec{x}_n(t)= \pmatrix x_n^1(t)\\ x_n^2(t) \endpmatrix=\pmatrix
x_0+\int_0^t a_n(x_n^1(\tau))d\tau+w(t)\\w(t)\endpmatrix
\mathop{\longrightarrow}\limits^W\vec{x}(\cdot)=\pmatrix
x_0+q\eta_\cdot+\widehat{w}(\cdot)\\ \widehat{w}(\cdot)
\endpmatrix
$$
when $n\to +\infty,$ $\{\widehat{w}(t)\}$ is a Wiener process.
\endproclaim

\demo{Proof} Without loss of generality we can assume that $x_0=0.$
Let us apply so-called "drift eliminating" transformation to the
first component of $\{\vec{x}_n(t)\}$ (see \cite {6}, p.111) :
$$
\vec{y}_n(t)=\pmatrix y_n^1(t) \\ y_n^2(t) \endpmatrix =\pmatrix
S_n(x_n^1(t))\\x^2_n(t)\endpmatrix=\pmatrix\int_0^t\sigma_n(y_n^1(
\tau))dw(\tau)\\ w(t)\endpmatrix,
$$
where
$$
S_n(x)=c\int_0^x\exp\{-2(A(nu)-A(0))\}du,\quad A(x)=\int_{-\infty}^x
a(z)dz,\quad c={\exp\{-2A(0)\}\over 1+\exp\{-2A\}},
$$
$\sigma_n(x)=S_n^\prime(S_n^{-1}(x))$ for all $x\in\Re.$ Then
$y^1_n (\cdot)\mathop{\longrightarrow}\limits^W y(\cdot),
n\to+\infty,$ where $\{y(t)\}$ is the solution of the following
SDE: $dy(t)=\sigma(y(t))d\widehat{w}(t),t\geq 0,\sigma(x)={1\over
2}(1-q\,\hbox{{\rm sign}}\, x),\{\widehat{w}(t)\}$ is a Wiener
process. The process $\{x^1(t)\},$ where $y(t)=r(x^1(t)),
r(x)={2x\over1-q\,\hbox{{\rm sign}}\, x}, x\in\Re$ is skew
Brownian motion.

The sequence $\{\vec{y}_n(t)\},n\geq 1$ is weakly compact because
each component of this sequence is weakly compact. Therefore we
prove the lemma if we show the uniqueness of the limit point. If
we prove that the equality $y(t)=\int_0^t\sigma(y(\tau))d\widehat{w}
(\tau)$ is valid for every limit point of $\{\vec{y}_n(t)\}$ then the
needed uniqueness follows from Nakao pathwise uniqueness theorem (see
\cite{8}). Note that $\sigma(\cdot)$ is separate from $0$ and has
bounded variation, i.e. Nakao theorem can be applied here.

Let $\{\vec{y}_{n_k} (t)\}$ be a convergent subsequence (we denote
it by $\{\vec{y}_k(t)\}$):
$$
\vec{y}_k(t)\mathop{\longrightarrow}\limits^W\vec{y}(t)\equiv\pmatrix
y(t)\\ \widehat{w}(t)
\endpmatrix,\quad k\to+\infty.
$$
Further we show that $y(t)=\int_0^t\sigma(y(\tau))d\widehat{w}
(\tau).$ For some $m\geq 1$ we denote by $\lambda_m$ the partition of
the segment $[0,t]: \lambda_m=\{0=t_0<t_1<\ldots <t_m=t\},$ where
$t_j=tj/m.$ Then
$$
{\Bbb E}\left|y(t)-\int_0^t\sigma(y(\tau))d\widehat{w}(\tau)\right|^2
\leq 2{\Bbb E}\left|y(t)-\sum_{j=0}^{m-1}\sigma(y(t_j))\Delta
\widehat{w}^j\right|^2+
$$
$$
+2{\Bbb E}\left|\int_0^t\sigma(y(\tau))d\widehat{w}(\tau)-
\sum_{j=0}^{m-1} \sigma(y(t_j))\Delta\widehat{w}^j\right|^2,\tag 1.1
$$
where $\Delta\widehat{w}^j=\widehat{w}(t_{j+1})-\widehat{w}(t_j).$

Consider the first summand. Since $\sigma(\cdot)$ has jump only at
one point and the process $\{y(t)\}$ has transition probability
density the mapping $\Phi:(y_k|_0^t,w|_0^t)\to y_k(t)-
\sum_{j=0}^{m-1}\sigma(y_k(t_j))\Delta w^j$ is continuous a.s. (we
denote by $\gamma|_0^t$ the trajectory of $\gamma$ on the segment
$[0,t]$). One can observe also that ${\Bbb E}\left|y_k(t)-
\sum_{j=0}^{m-1}\sigma(y_k(t_j))\Delta w^j\right|^4\leq 48t^2.$
Therefore using the theorems 5.1 and 5.4, \cite{9} we get the
equality:
$$
{\Bbb E}\left|y(t)-\sum_{j=0}^{m-1}\sigma(y(t_j))\Delta\widehat{w}^j
\right|^2=\lim_{k\to+\infty}{\Bbb E}\left|y_k(t)-\sum_{j=0}^{m-1}
\sigma(y_k(t_j))\Delta w^j\right|^2.\tag 1.2
$$
Let us denote $\widehat{y}_k(t)=y_k(t_j),t\in [t_j,t_{j+1}).$ We
obtain
$$
{\Bbb E}\left|y_k(t)-\sum_{j=0}^{m-1}\sigma(y_k(t_j))\Delta w^j
\right|^2={\Bbb E}\left|\int_0^t\left[\sigma_k(y_k(\tau))-\sigma
(\widehat{y}_k(\tau))\right]dw(\tau)\right|^2=
$$
$$
={\Bbb E}\int_0^t\left[\sigma_k(y_k(\tau))-\sigma(\widehat{y}_k
(\tau))\right]^2 d\tau\leq 2{\Bbb E}\int_0^t\left[\sigma_k(y_k(\tau))
-\sigma(y_k(\tau))\right]^2d\tau+
$$
$$
+2{\Bbb E}\int_0^t\left[\sigma(y_k(\tau))-\sigma(\widehat{y}_k(\tau))
\right]^2d\tau.\tag 1.3
$$

Let us estimate the first summand in (1.3). Put $\tau_R=\inf\{t\geq
0:|y_k(t)|\geq R\},$ one has
$$
{\Bbb E}\int_0^t\left[\sigma_k(y_k(\tau))-\sigma(y_k(\tau))\right]^2
d\tau={\Bbb E} {\hbox{\rm 1\!I}}_{\{\tau_R\leq t\}} \int_0^t\left[
\sigma_k(y_k(\tau))-\sigma(y_k(\tau))\right]^2d\tau+
$$
$$
+{\Bbb E}{\hbox{\rm 1\!I}}_{\{\tau_R>t\}}\int_0^t\left[\sigma_k
(y_k(\tau))-\sigma(y_k(\tau))\right]^2 d\tau\leq t\max_{x,y\in
\Re}[\sigma_k(x)-\sigma(y)]^2{\Bbb P}\{\tau_R\leq t\}+
$$
$$
+{\Bbb E}\int_0^{t\wedge\tau_R}\left[\sigma_k(y_k(\tau))-\sigma
(y_k(\tau))\right]^2 d\tau.\tag 1.4
$$
For some $\delta>0$ take $R=R_\delta=\sqrt{4t\over \delta}$ such
that the following inequality holds
$$
{\Bbb P}\{\tau_{R_\delta}<t\}={\Bbb P}\{\max_{0\leq s\leq t}
|y_k(t)|\geq R_\delta\}\leq {1\over R_\delta^2}{\Bbb E} \left[
\max_{0\leq s\leq t}|y_k(t)|\right]^2\leq
$$
$$
\leq {4\over R^2_\delta}{\Bbb E}\int_0^t\sigma_k^2(y_k(s))ds\leq
{4t\over R_\delta^2}=\delta.\tag 1.5
$$
It follows from Krylov's inequality (see, for example, \cite{10},
lemma 1, p.562) that there exists a constant $q_{\delta,t}$ such
that the following estimate holds
$$
{\Bbb E}\int_0^{t\wedge\tau_{R_\delta}} \left[\sigma_k(y_k(\tau))-
\sigma(y_k(\tau))\right]^2d \tau\leq q_{\delta,t}\left[\int_0^t
d\tau\int_{|x|\leq R_\delta}[\sigma_k(x)-\sigma(x)]^4dx\right]^{1/2}.
$$

Let $k\to+\infty.$ Using the Lebesgue's majorized convergence theorem
we see that the second summand in (1.4) tends to $0.$ Therefore, we
get
$$
\limsup_{k\to+\infty}{\Bbb E}\int_0^t\left[\sigma_k(y_k(\tau))-
\sigma(y_k(\tau))\right]^2d\tau \leq 4t\delta.
$$

Then we proceed to the limit as $\delta\to 0$ and obtain that the
first summand in (1.3) tends to $0$ when $k\to +\infty.$

Consider the second summand on the right hand side of (1.3). Using
the explicit form of the function $\sigma(\cdot)$ we get
$$
{\Bbb E}\int_0^t\left[\sigma(y_k(\tau))-\sigma(\widehat{y}_k(\tau))
\right]^2d\tau={q^2\over 4}{\Bbb E}\int_0^t\left[\,\hbox{{\rm
sign}}\, y_k(s)-\,\hbox{{\rm sign}}\,\widehat{y}_k(s)\right]^2ds=
$$
$$
={q^2\over 4}\sum_{j=0}^{m-1}{\Bbb E} \int_{t_j}^{t_{j+1}}
\left[\,\hbox{{\rm sign}}\,y_k(s)-\,\hbox{{\rm sign}}\,
y_k(t_j)\right]^2ds=
$$
$$
=q^2\sum_{j=0}^{m-1}\int_{t_j}^{t_{j+1}}{\Bbb P}\left\{\,\hbox{{\rm
sign}}\, y_k(s)\not=\,\hbox{{\rm sign}}\, y_k(t_j)\right\}ds\leq
q^2\sum_{j=0}^{m-1}\int_{t_j}^{t_{j+1}}{\Bbb P}\left\{|y_k(s)|\leq
\varepsilon\right\}ds+
$$
$$
+q^2\sum_{j=0}^{m-1}\int_{t_j}^{t_{j+1}}{\Bbb P}\left\{\,\hbox{{\rm
sign}}\, y_k(s)\not=\,\hbox{{\rm sign}}\, y_k(t_j),|y_k(s)|>
\varepsilon\right\}ds, \tag 1.6
$$
here some $\varepsilon>0$ is fixed. Observing that $\{\,\hbox{{\rm
sign}}\, y_k(s)\not=\,\hbox{{\rm sign}}\, y_k(t_j),$ $|y_k(s)|>
\varepsilon\}\subseteq\{\max_{t_j\leq s\leq t_{j+1}} |y_k(s)-
y_k(t_j)|>\varepsilon\}$ we estimate the second summand in (1.6) in
the following way:
$$
{\Bbb E}\int_{t_j}^{t_{j+1}}\!\!\!\!{\hbox{\rm 1\!I}}_{\{\,\hbox{{\rm
sign}}\, y_k(s)\not=\,\hbox{{\rm sign}}\,y_k(t_j), |y_k(s)|>
\varepsilon\}}ds\leq {t\over m}{\Bbb P}\{\max_{t_j\leq s\leq t_{j+1}}
|y_k(s)-y_k(t_j)|>\varepsilon\}\leq
$$
$$
\leq {t\over m\varepsilon^2}{\Bbb E}\left[\max_{t_j\leq s\leq
t_{j+1}}|y_k(s)- y_k(t_j)|\right]^2\leq{4t\over m\varepsilon^2}
{\Bbb E}\int_{t_j}^{t_{j+1}}\sigma_k^2(y_k(s))ds\leq{4t^2\over
m^2\varepsilon^2}.\tag 1.7
$$

Consider the first summand on the right hand side of (1.6). We
have
$$
\sum_{j=0}^{m-1}\int_{t_j}^{t_{j+1}}{\Bbb P}\left\{|y_k(s)|
\leq\varepsilon\right\}ds =\int_0^t{\Bbb P}\left\{|y_k(s)|\leq
\varepsilon\right\}=\int_0^t{\Bbb P}\left\{|y_k(s)| \leq\varepsilon,
\tau_R<t\right\}ds+
$$
$$
+\int_0^t{\Bbb P}\left\{|y_k(s)|\leq\varepsilon,\tau_R\geq t
\right\}ds\leq t{\Bbb P}\{ \tau_R<t\}+{\Bbb E}\int_0^{t\wedge\tau_R}
{\hbox{\rm 1\!I}}_{\{|y_k(s)|\leq\varepsilon\}}ds.\tag 1.8
$$

Let $R_\delta$ be as in previous one. It follows from Krylov's
inequality that there exists a constant $q_{\delta,t}$ such that
the following estimate holds
$$
{\Bbb E}\int_0^{t\wedge\tau_{R_\delta}}{\hbox{\rm 1\!I}}_{\{
|y_k(s)|\leq\varepsilon\}}ds\leq q_{t,\delta} \left[\int_0^{t}
ds\int_{|x|<R_\delta}{\hbox{\rm 1\!I}}_{\{|x|\leq\varepsilon\}}
dx\right]^{1/2}=2q_{t,\delta}\sqrt{t(\varepsilon\wedge R_\delta)}.
\tag 1.9
$$
For the second summand in (1.1) one can write estimates analogous to
(1.6)-(1.9) and obtain inequality
$$
{\Bbb E}\int_0^t\left[\sigma(y(\tau))-\sigma(\widehat{y}(\tau))
\right]^2 d\tau\leq {4t^2q^2\over m\varepsilon^2}+t\delta
+2q_{t,\delta}\sqrt{t(\varepsilon\wedge R_\delta)},\tag 1.10
$$
with the same $\varepsilon,\delta,q_{t,\delta},R_\delta.\quad
\widehat{y}(t)=y(t_j),t\in[t_j, t_{j+1}).$

Finally, from (1.1), (1.5)-(1.10) we get
$$
{\Bbb E}\left|y(t)-\int_0^t\sigma(y(\tau))d\tau\right|^2\leq
{24t^2q^2\over m\varepsilon^2}+6t\delta+12q_{t,\delta}
\sqrt{t(\varepsilon\wedge R_\delta)}.
$$
Proceeding first $m\to+\infty,$ then $\varepsilon\to 0$ and, at
last, $\delta\to 0,$ we obtain the required result. The lemma is
proved.\enddemo

\proclaim{Lemma 2} Consider the sequence of processes in $\Re^2:$
$$
\vec{x}_n(t)=\pmatrix x_n^1(t) \\ x_n^2(t)\endpmatrix =\pmatrix x^1_0
+\int_0^t a^1_n(x_n^1(\tau))d\tau+w(t) \\ x^2_0+\int_0^t a^2_n
(x_n^2(\tau))d\tau+w(t)\endpmatrix,
$$
where each component $\{x^i_n(t)\},n\geq 1,i=1,2$ are defined as in
Proposition 1. Then
$$
\vec{x}_n(\cdot)\mathop{\longrightarrow}\limits^W\vec{x}(\cdot)=
\pmatrix x^1(\cdot)\\x^2(\cdot)\endpmatrix=\pmatrix
x^1_0+q^1\eta_\cdot^1+\widehat{w}(\cdot)\\ x^2_0+
q^2\eta_\cdot^2+\widehat{w}(\cdot)\endpmatrix,\quad n\to+\infty,
$$
$\{x^i(t)\},i=1,2$ are skew Brownian motions constructed as the
functional of the same Wiener process $\{\widehat{w}(t)\}$ and
$q^1, \{\eta_t^1\},q^2, \{\eta_t^2\}$ defined as in Proposition 1.
\endproclaim

\demo{Proof}The sequence of the processes $\vec{x}_n(\cdot)$ is
weakly compact. This mean that the sequence $\vec{x}_n(\cdot)$ has a
limit point. If we show the uniqueness of this point then we prove
this lemma. We prove the uniqueness by contradiction. Let
$\vec{x}_{n_k^1}(\cdot)$ and $\vec{x}_{n_k^2}(\cdot)$ be converged
subsequences of $\vec{x}_n(\cdot)$ with different limit points.
Consider sequences in $\Re^3:$
$$
\vec{X}^1_k(\cdot)=\pmatrix x_{n_k^1}^1(\cdot)\\x_{n_k^1}^2(\cdot)\\
\widehat{w}(\cdot)\endpmatrix,\quad \vec{X}^2_k(\cdot)=\pmatrix
x_{n_k^2}^1(\cdot)\\x_{n_k^2}^2(\cdot)\\\widehat{w}(\cdot)\endpmatrix
.
$$
Analogously to previous consideration these sequences are weakly
compact. Let $\{\vec{X}^1_{m_k^1}(\cdot)\}$ and $\{\vec{X}^2_{m_k^2}
(\cdot)\}$ be some convergent subsequences of the sequences
$\{\vec{X}^1_k(\cdot)\}$ and $\{\vec{X}^2_k(\cdot)\}:$
$$
\vec{X}^i_{m_k^i}(\cdot)=\pmatrix x_{m_k^i}^1(\cdot)\\x_{m_k^i}^2
(\cdot)\\\widehat{w}(\cdot)\endpmatrix \mathop{\longrightarrow}
\limits^W\vec{X}_i(\cdot)=\pmatrix x^1_i(\cdot)\\x_i^2(\cdot)
\\x_i^3(\cdot)\endpmatrix, \quad i=1,2.
$$

Consider the process $\{\vec{X}_1(\cdot)\}.$ According to Lemma 1
the first component $x_1^1(\cdot)$ is skew Brownian motion
constructed as the functional of the Wiener process $x_1^3(\cdot),$
i.e. there exists a measurable functional $\Phi_{q^1,x_0^1}:
C[0,+\infty)\to C[0,+\infty)$ such that $x_1^1(\cdot)=
\Phi_{q^1,x_0^1}(x_1^3(\cdot)).$ The second component is the same
functional of the $x_1^3(\cdot): x_1^2(\cdot)=\Phi_{q^2,x_0^2}
(x_1^3(\cdot)).$ Therefore the distribution of the process
$\{\vec{X}_1(\cdot)\}$ is the image of the Wiener measure under the
mapping $\Psi:C[0,+\infty)\to C([0,+\infty),\Re^3),$ where
$$
\Psi:y(\cdot)\to\pmatrix\Phi_{q^1,x_0^1}(y(\cdot))\\\Phi_{q^2,x_0^2}(
y(\cdot))\\y(\cdot)\endpmatrix.
$$

The same arguments are valid for the process $\{\vec{X}_2(\cdot)\}.$
This means that the distributions of the processes $\{\vec{X}_1(
\cdot)\}$ and $\{\vec{X}_2(\cdot)\}$ coincide, that gives
contradiction. The lemma is proved.

\enddemo

\demo{Proof of Theorem 1} Let $a_1(\cdot)$ be a function satisfying
the conditions of Proposition 1 and let $a^2(x)=a^1(x)+{A_2-A_1
\over\sqrt{2\pi}}\exp\{-{x^2\over 2}\},x\in\Re,$ $\,\hbox{{\rm th}}\,
A_i=q_i, i=1,2.$ One can see that $a^2(\cdot)$ satisfies the
conditions of Proposition 1 too and $a_n^2(\cdot)\geq a_n^1(\cdot)$
for all $n\geq 1,$ where $a_n^i(\cdot)=na_i(n\cdot).$ For
$a_n^1(\cdot),a_n^2 (\cdot)$ consider the sequence of diffusion
processes $\{\vec{x}_n (\cdot)\}=\{(x_n^1 (\cdot),x_n^2(\cdot))\}$
defined as in Lemma 2. These processes satisfies the conditions of
the comparison theorem for diffusion processes (see, for example,
\cite {3}, Section VI, Theorem 1.1), i.e. ${\Bbb P}\{x^n_1(t)\leq
x^n_2(t),\forall t\geq 0\}=1,\quad n\geq 1.$

The set $A=\left\{(\varphi_1(\cdot),\varphi_2(\cdot))\in C([0,
+\infty),\Re^2):\ \varphi_1(t) \leq \varphi_2(t),\forall t\geq
0\right\}$ is closed in \ $C([0,+\infty),\Re^2).$ Therefore from the
properties of weak convergence we have ${\Bbb P}\{x_1(t)\leq x_2(t),$
$\forall t\geq 0\}\geq \limsup_{n\to+\infty}{\Bbb P}\{x^1_n(t)\leq
x^2_n(t),\forall t\geq 0\}=1.$ The theorem is proved. \enddemo

\proclaim{Corollary 1} Consider a pair of skew Brownian motions
$\{x^1(t)\},\{x^2(t)\}$ constructed as the functional of the Wiener
process $\{w(t)\}$ with different skewing parameters $q_1,q_2
\in(-1,1)$ and started from the same point $x\in\Re.$ Then the
equality
$$
{\Bbb E}|x^1(t)-x^2(t)|=|q_1-q_2|I_t(x),\quad t\geq 0,
$$
holds with $I_t(x)={\Bbb E}_x\eta_t^1={\Bbb E}_x\eta_t^2=\int_0^t
{1\over\sqrt{2\pi \tau}}\exp\left\{-{x^2\over 2\tau}\right\}d\tau.$
\endproclaim

\proclaim{Corollary 2} Consider a pair of skew Brownian motions
$\{x^1(t)\},\{x^2(t)\}$ constructed as the functional of the Wiener
process $\{w(t)\}$ with the same skewing parameters $q\in(-1,0)\cup
(0,1)$ and started from the different points $x_0^1,x_0^2\in\Re.$
Then the inequalities
$$
{\Bbb E}|x^1(t)-x^2(t)|\leq|x^1_0-x_0^2|+|q||I_t(x_0^1)-I_t(x_0^2)|,
$$
$$
{\Bbb E}|\eta_t^1(t)-\eta_t^2(t)|\leq {1\over |q|}|x_0^1-x_0^2|+
|I_t(x_0^1)-I_t(x_0^2)|, \quad t\geq 0,\tag 1.11
$$
hold, where the function $I_t(\cdot), t\geq 0$ is defined in the
Corollary 1.\endproclaim

The proofs of Corollary 1 and Corollary 2 are easy and omitted.

\remark{Remark 1} Consider the case of $|q|=1.$ We assume that the
phase space is equal $[0,+\infty)$ when $q=1$ and $(-\infty,0]$
when $q=-1.$ The estimate similar to (1.11) in this situation also
holds true (see \cite{11}):
$$
{\Bbb E}|x_1(t)-x_2(t)|^2\leq|x_0^1-x_0^2|^2.
$$\endremark

\remark{Remark 2} In the case of $q=0$ the inequality
$$
{\Bbb E}|\eta_t^1-\eta_t^2|^2\leq 16|x_0^1-x_0^2|^2+{8\sqrt{t}\over
\sqrt{\pi}}|x_0^1-x_0^2|
$$
holds. The proof of this remark is easy corollary of Tanaka's formula
for local time of Wiener process.\endremark

\head 2. On stochastic continuity of strong solution to SDE with
singular coefficients\endhead

Let $S$ be a hyperplane in $\Re^d$ orthogonal to the fixed ort
$\nu\in\Re^d.$ We denote by $\pi_S$ the operator of orthogonal
projection on $S.$ For a pair of independent Wiener processes
$\{w(t)\}$ and $\{\widetilde{w}(t)\}$ in $\Re^d$ and $S$
respectively, for given parameters $x_0\in\Re^d, q\in[-1,1],$
given measurable function $\alpha:S\to S$ and operator $\beta:S\to
{\Cal L}_+(S)$ (${\Cal L}_+(S)$ is the space of all linear symmetric
nonnegative operators on $S$) we consider the following stochastic
equation in $\Re^d$
$$
x(t)=x_0+\int_0^t{\left(q\nu+\alpha(x^S(\tau))\right)d\eta_\tau}+
\int_0^t{\widetilde{\beta}(x^S(\tau))d\widetilde{w}(\eta_\tau)}+w(t)
\tag 2.1
$$
where $\widetilde{\beta}(\cdot)=\beta^{1/2}(\cdot), x^S(\cdot)=\pi_S
x^S(\cdot).$ It is proved in \cite {12} that under the following
assumptions on the coefficients
\roster
\item
$\sup_{x\in S}(|\alpha(x)|+\left\|\widetilde{\beta}(x)\right\|)\leq
K$,
\item
$|\alpha(x)-\alpha(y)|^2+\left\|\widetilde{\beta}(x)-\widetilde{\beta
}(y)\right\|^2\leq K|x-y|^2,$ for all $x,y\in S$
\endroster
for some $K>0$ the solution to the equation (2.1) exists and is
unique. In the next theorem we prove that this solution continuously
depends on the starting point.

\proclaim {Theorem 2} Let $\{x_n(t)\},n\geq 1$ be the sequence of
the solutions to (2.1) started from $\{x_n\}\subset\Re^d$ and let
$x_n\to x\in\Re^d$ when $n\to+\infty.$ Then for all $t\geq 0$
$x_n(t)\mathop{\longrightarrow}\limits^P x(t), n\to+\infty$ where
$\{x(t)\}$ is the solution to (2.1) started from $x.$\endproclaim

\demo{Proof} Consider a new process $\rho_t=\inf\{s\geq 0:\eta_s\geq
t\}.$ It is nonnegative left continuous increasing process,
$\rho_0=0$. Also $\rho_t\to +\infty$ when $t\to +\infty$ and
$\eta_{\rho_t}=t$ for all $t\geq 0.$ Let us note that $\{\rho_t\}$ is
the stopping time w.r.t. ${\Cal F}_t^w=\sigma\{w(u),0\leq u\leq t\}.$

We substitute the process $\rho_t$ instead of $t$ in equation (2.1)
and set $u=\eta_\tau$ in all integrals in (2.1). Then we obtain
$$
x(\rho_t)=x+\int_0^t\left(q\nu+\alpha(x^S(\rho_u))\right)du+\int_0^t
{\widetilde{\beta}(x^S(\rho_u))d\widetilde{w}(u)}+w(\rho_t)\tag 2.2
$$
We construct the processes $\{\rho_t^n\},\{x_n(\rho_t^n)\}$ in the
same way.

\proclaim {Lemma 3}For all $t\geq 0:$ $\rho^n_t\mathop
{\longrightarrow}\limits^P\rho_t$ when $n\to+\infty.$\endproclaim

\demo{Proof} Firstly we prove that ${\Bbb P}\{\rho_t^n-\rho_t>
\varepsilon\}\to 0,n\to+\infty.$ We can write
$$
{\Bbb P}\{\rho_t^n-\rho_t>\varepsilon\}={\Bbb P}\{\rho_t^n>
\varepsilon+\rho_t\}=\sum_{k\geq 0}{\Bbb P}\left\{\rho_t^n>
\varepsilon+\rho_t,\rho_t^n\in\left(\left.{(k+1)\varepsilon
\over 2},{(k+2)\varepsilon \over 2}\right]\right.\right\}\leq
$$
$$
\leq\sum_{k\geq 0}{\Bbb P}\left\{\rho_t<{k\varepsilon\over 2},
\rho_t^n>{(k+1)\varepsilon\over 2}\right\}= \sum_{k\geq 0}{\Bbb P}
\left\{\eta_{k\varepsilon/2}>t,\eta^n_{(k+1)\varepsilon/2}<t\right\}.
$$

Consider the k-th summand. For some $\delta>0$
$$
{\Bbb P}\left\{\eta_{k\varepsilon/2}>t,\eta^n_{(k+1)\varepsilon/2}<t
\right\}={\Bbb P}\left\{\eta_{k\varepsilon/2}>t+\delta,
\eta^n_{(k+1)\varepsilon/2}<t\right\}+
$$
$$
+{\Bbb P}\left\{t<\eta_{k\varepsilon/2}\leq t+\delta,
\eta^n_{(k+1)\varepsilon/2}<t\right\} \leq {\Bbb P}\left\{\eta_{(k+1)
\varepsilon/2}>t+\delta,\eta^n_{(k+1)\varepsilon/2}<t\right\}+
$$
$$
+{\Bbb P}\left\{t<\eta_{k\varepsilon/2}\leq t+\delta\right\}\leq
{\Bbb P}\left\{|\eta_{(k+1)\varepsilon/2}-\eta^n_{(k+1)
\varepsilon/2}|>\delta\right\}+{\Bbb P}\left\{t<\eta_{k\varepsilon/2}
\leq t+\delta\right\}. \tag 2.3
$$

We use the distribution of $\{\eta_t\}$ for estimating the second
summand in (2.3). We remind that $\{\eta_t\}$ has the same
distribution with the local time in $0$ of a Wiener process in
$\Re$ started from $x^\nu=(x,\nu)$:
$$
{\Bbb P}\{\eta_t^x<a\}= \left(1-2\int_{{|x^\nu|+a\over
\sqrt{t}}}^{+\infty}{\exp\left\{-u^2/2\right\}\over\sqrt{2\pi}}du
\right){\hbox{\rm 1\!I}}_{\{a>0\}}.
$$

Thus we have
$$
{\Bbb P}\left\{t<\eta_{k\varepsilon/2}\leq t+\delta\right\}=
2\int_{{|x^\nu|+t\over\sqrt{k\varepsilon/2}}}^{{|x^\nu|+(t+\delta)
\over\sqrt{k\varepsilon/2}}}{\exp\left\{-u^2/2\right\}\over\sqrt{2\pi
}}du\leq 2\sqrt{{2\over\pi k\varepsilon}}\delta.\tag 2.4
$$

The first summand in (2.3) is estimated by using Chebyshev's
inequality and Corollary 2 in the case of $q\in(-1,0)\cup(0,1):$
$$
{\Bbb P}\left\{|\eta_{(k+1)\varepsilon/2}-\eta^n_{(k+1)
\varepsilon/2}|>\delta\right\}\leq {1\over \delta} {\Bbb E}
|\eta_{(k+1)\varepsilon/2}-\eta^n_{(k+1)\varepsilon/2}|\leq
$$
$$
\leq {1\over \delta}\left(c_{1,q}|x_n-x|+ c_{2,q}|I_{(k+1)
\varepsilon/2}(x_n)-I_{(k+1)\varepsilon/2}(x)|\right). \tag 2.5
$$

Finally we obtain from (2.3)-(2.5) that
$$
{\Bbb P}\left\{\rho_t^n>\varepsilon+\rho_t,\rho_t^n\in
\left(\left.{(k+1)\varepsilon\over 2},{(k+2) \varepsilon\over
2}\right]\right.\right\}\leq\left(2\sqrt{{2\over\pi k\varepsilon}}
\delta+\right.
$$
$$
\left.+{1\over\delta}c_{1,q}|x_n-x|+c_{2,q}|I_{(k+1)\varepsilon/2}
(x_n)-I_{(k+1)\varepsilon/2}(x)|\right).
$$

Proceeding first $n\to +\infty,$ then $\delta\to 0$ we see that
for all $k\geq 0:$
$$
{\Bbb P}\left\{\rho_t^n>\varepsilon+\rho_t,\rho_t^n\in\left(
\left.{(k+1)\varepsilon\over2},{(k+2)\varepsilon \over
2}\right]\right.\right\}\to 0, \quad n\to+\infty.
$$

Taking into account that
$$
\sum_{k\geq 0}{\Bbb P}\left\{\rho_t^n>\varepsilon+\rho_t,\rho_t^n\in
\left(\left.{(k+1)\varepsilon \over 2},{(k+2)\varepsilon \over
2}\right]\right.\right\}\leq
$$
$$
\leq \sum_{k\geq 0}{\Bbb P}\left\{\rho_t^n\in\left(\left.
{(k+1)\varepsilon \over 2},{(k+2)\varepsilon \over 2}\right]
\right.\right\}\leq {\Bbb P}\left\{\rho_t^n\geq 0\right\}=1
$$
we see that conditions of Lebesgue's majorized convergence theorem
is satisfied. Therefore  ${\Bbb P}\{\rho^n_t-\rho_t>\varepsilon\}\to
0, n\to+\infty.$ The same arguments can be made in the case $|q|=1$
(by using Remark 1), $q=0$ (by using Remark 2). We prove that ${\Bbb
P}\{\rho_t-\rho^n_t>\varepsilon\}\to 0, n\to+\infty$ in the same
way. The lemma is proved.

\proclaim{Lemma 4} For all $t\geq 0:$ $|x_n(\rho_t^n)-x(\rho_t)|
\mathop{\longrightarrow}\limits^P 0, n\to+\infty.$ \endproclaim

\remark{Remark 3} One can see that ${\Bbb E}|x_n(\rho_t^n)-
x(\rho_t)|^p= +\infty, p\geq 1$ because ${\Bbb E}\rho_t=+\infty.$
Thus we cannot apply here standard technique such as using martingale
inequalities.\endremark

\demo{Proof} For given $N>0, C>0, n\geq 1$ we consider the random set
$$
A_{N,C}^{n,t}=\left\{\omega\in\Omega: \rho_t<N,\rho_t^n<N,
\left\|w\right\|_{\,\hbox{{\rm Hol}}\,_{1/4}([0,N])}\leq C\right\},
$$
where $\left\|\cdot\right\|_{\,\hbox{{\rm Hol}}\,_{1/4}([0,N])}$
is the H\"older norm with parameter $1/4$ on the segment $[0,N].$
Note that $A_{N,C}^{n,t} \in{\Cal F}^w_N$ and $A_{N,C}^{n,t}\subseteq
A_{N,C}^{n,s}$when $s\leq t.$

Then for all $\varepsilon>0, t\in[0,T]$ we get
$$
{\Bbb P}\{|x_n(\rho^n_t)-x(\rho_t)|>\varepsilon\}\leq{\Bbb P}
\{\{|x_n(\rho^n_t)-x(\rho_t)|>\varepsilon\}\bigcap
A_{N,C}^{n,t}\}+{\Bbb P}\{\rho_t\geq N\}+
$$
$$
+{\Bbb P}\{\rho_t^n\geq N\}+{\Bbb P}\{\left\|w\right\|_{\,\hbox{{\rm
Hol}}\,_{1/4} ([0,N])}>C\}\tag 2.6
$$

Let us estimate the second moment of the process $\{(x_n(\rho^n_t)-
x(\rho_t)){\hbox{\rm 1\!I}}_{A_{N,C}^{n,t}}\}$:
$$
{\Bbb E}|x_n(\rho^n_t)-x(\rho_t)|^2{\hbox{\rm 1\!I}}_{A_{N,C}^{n,t}}
\leq 4|x_n-x|^2+ 4t\int_0^t{{\Bbb E}|\alpha(x_n^S(\rho^n_u))-
\alpha(x^S(\rho_u))|^2 {\hbox{\rm 1\!I}}_{A_{N,C}^{n,t}}du}+
$$
$$
+4{\Bbb E}\left[\int_0^t{(\widetilde{\beta}(x_n^S(\rho^n_u))-
\widetilde{\beta}(x^S(\rho_u)))d\widetilde{w}(u)}\right]^2{\hbox{\rm
1\!I}}_{A_{N,C}^{n,t}}+4{\Bbb E}|w(\rho^n_t)-w(\rho_t)|^2{\hbox{\rm
1\!I}}_{A_{N,C}^{n,t }}\leq
$$
$$
\leq 4|x_n-x|^2+4KT{\Bbb E}\int_0^t|x_n^S(\rho^n_u)-x^S(\rho_u)|^2
{\hbox{\rm 1\!I}}_{A_{N,C}^{n,u}}du+
$$
$$
+4{\Bbb E}\left[\int_0^t{(\widetilde{\beta}(x_n^S(\rho^n_u))-
\widetilde{\beta}(x^S(\rho_u)))d\widetilde{w}(u)}{\hbox{\rm
1\!I}}_{A_{N,C}^{n,t}}\right]^2+4C{\Bbb E}|\rho^n_t-\rho_t|^{1/2}
{\hbox{\rm 1\!I}}_{A_{N,C}^{n,t}}.\tag 2.7
$$

Consider the third summand. Let us put $\widehat{{\Cal F}}_t={\Cal
F}^w_\infty\vee {\Cal F}^{\widetilde{w}}_t.$ One can observe that
the process $\{\widetilde{\beta}(x_n^S (\rho^n_t))- \widetilde{\beta}
(x^S(\rho_t))\}$ is $\{\widehat{{\Cal F}}_t\}-$adapted, the process
$\{\widetilde{w}(t)\}$ is the Wiener process w.r.t. $\{\widehat{{\Cal
F}}_t\}$ and $A_{N,C}^{n,t}\in\widehat{{\Cal F}}_t$ for all $t\geq
0.$ Therefore the equality
$$
{\Bbb E}\int_0^t(\widetilde{\beta}(x_n^S(\rho^n_u))-
\widetilde {\beta}(x^S(\rho_u)))d\widetilde{w}(u){\hbox{\rm
1\!I}}_{A_{N,C}^{n,t}}={\Bbb E}\int_0^t(\widetilde{\beta}
(x_n^S(\rho^n_u))-\widetilde{\beta}(x^S(\rho_u)))
{\hbox{\rm 1\!I}}_{A_{N,C}^{n,t}}d\widetilde{w}(u)
$$
holds. Thus we obtain
$$
{\Bbb E}\left[\int_0^t{(\widetilde{\beta}(x_n^S(\rho^n_u))-
\widetilde{\beta}(x^S(\rho_u))){\hbox{\rm
1\!I}}_{A_{N,C}^{n,t}}d\widetilde{w}(u)}\right]^2=
$$
$$
=\int_0^t{\Bbb E}\left\|\widetilde{\beta}(x_n^S(\rho^n_u))
-\widetilde{\beta}(x^S(\rho_u))\right\|^2 {\hbox{\rm 1\!I}}_{A_{N,C
}^{n,t}}du\leq K\int_0^t {\Bbb E}|x_n^S(\rho^n_u)-x^S(\rho_u)|^2
{\hbox{\rm 1\!I}}_{A_{N,C}^{n,u}}du. \tag 2.8
$$

It follows from (2.7) and (2.8) that
$$
{\Bbb E}|x_n(\rho^n_t)-x(\rho_t)|^2{\hbox{\rm 1\!I}}_{A_{N,C}^{n,t}}
\leq 4|x_n-x|^2+
$$
$$
+4K(T+1){\Bbb E}\int_{0}^{t}{|x_n^S(\rho^n_u)-x^S(\rho_u)
|^2{\hbox{\rm 1\!I}}_{A_{N,C}^{n,u}} du}+4C{\Bbb E}|\rho^n_t-
\rho_t|^{1/2}{\hbox{\rm 1\!I}}_{A_{N,C}^{n,t}}.
$$

Using the Grownwall-Bellman inequality we obtain
$$
{\Bbb E}|x_n(\rho^n_t)-x(\rho_t)|^2{\hbox{\rm
1\!I}}_{A_{N,C}^{n,t}}\leq 4|x_n-x|^2+ 4C{\Bbb
E}|\rho^n_t-\rho_t|^{1/2}{\hbox{\rm 1\!I}}_{A_{N,C}^{n,t}}+
$$
$$
+4K(T+1)\int_0^t\exp\left\{4K(T+1)(t-u)\right\}\left(4|x_n-x|^2+
4C{\Bbb E}|\rho^n_u-\rho_u|^{1/2}{\hbox{\rm
1\!I}}_{A_{N,C}^{n,u}}\right)du \tag 2.9
$$
for all $t\in[0,T].$

It follows from Lemma 3 and from the fact that the processes
$\{\rho_t^n{\hbox{\rm 1\!I}}_{\{A_{N,C}^{n,t}\}}\},\{\rho_t{\hbox{\rm
1\!I}}_{\{A_{N,C}^{n,t}\}}\}$ are bounded for fixed $N,C,t$ that
${\Bbb E}|\rho_t^n-\rho_t|^{1/2}{\hbox{\rm 1\!I}}_{A_{N,C}^{n,t}}\to
0, n\to+\infty.$ From (2.9) we see that for all $\delta>0$ there
exists $n_0>0$ such that ${\Bbb E}|x_n(\rho_t^n)-x(\rho_t)|^2
{\hbox{\rm 1\!I}}_{A_{N,C}^{n,t}}<\delta$ for all $n\geq n_0.$ Using
the Chebyshev's inequality for the first summand of (2.6) we obtain
that
$$
{\Bbb P}\{|x_n(\rho^n_t)-x(\rho_t)|>\varepsilon\}\leq{\delta\over
\varepsilon^2}+{\Bbb P}\{\rho_t\geq N\}+\sup_{t\in[0,T]}{\Bbb P}
\{\rho_t^n\geq N\}+{\Bbb P}\{\left\|w\right\|_{ \,\hbox{{\rm
Hol}}\,_{1/4}([0,N])}>C\}
$$
for all $n\geq n_0.$ Proceeding first $\delta\to 0,$ then $C\to
+\infty$ and, at last, $N\to+\infty,$ we obtain that ${\Bbb P}
\{|x_n(\rho^n_t)- x(\rho_t)|>\varepsilon\}\to 0, n\to+\infty.$
Note that ${\Bbb P}\{\rho_t\geq N\}\leq {2(|x^\nu|+T)\over
\sqrt{2\pi N}}\to 0, N\to+\infty.$ Lemma is proved.

Let us return to the proof of the Theorem 2. We have
$$
{\Bbb P}\{|x_n(t)-x(t)|>\varepsilon\}\leq{\Bbb P} \{|x_n(t)-x(t)-
x_n(\rho^n_{\eta_t^n})+x(\rho_{\eta_t})|>\varepsilon/2\}+
$$
$$
+{\Bbb P}\{|x_n(\rho^n_{\eta_t^n})-x(\rho_{\eta_t})|>\varepsilon/2\}
\tag 2.10
$$
It follows from definition of $\{\rho_t\}$ that $\rho_{\eta_t}\leq
t.$ One can observe that $\eta_s=const, s\in[\rho_{\eta_t},t],$ thus
$x(s)=x_0+ w(s), s\in[\rho_{\eta_t},t]$ and $|x_n(t)-x(t)-x_n(
\rho^n_{\eta_t^n})+x(\rho_{\eta_t})|=|w(\rho^n_{\eta_t^n})-
w(\rho_{\eta_t})|.$

Consider the second summand in (2.10). We can write
$$
|x_n(\rho^n_{\eta_t^n}\!)-x(\rho_{\eta_t}\!)|\leq |x_n-x|+\left|
\int_{\eta_t}^{\eta_t^n}\!\!\!\!\!\left(q\nu+\alpha(x_n^S(\rho_u^n))
\right)du\right|+\left|\int_0^{\eta_t}\!\!\!\!\!\left(\alpha(x_n^S(
\rho_u^n)-\alpha(x^S(\rho_u))\right)du\right|+
$$
$$
+\left|\int_{\eta_t}^{\eta_t^n}\widetilde{\beta}(x_n^S(\rho_u^n))
d\widetilde{w}(u)\right|+\left|\int_0^{\eta_t}\left(\widetilde{\beta}
(x_n^S(\rho_u^n))-\widetilde{\beta}(x^S(\rho_u))\right)d\widetilde{w}
(u)\right|+\left|w(\rho^n_{\eta_t^n})-w(\rho_{\eta_t})\right|.
\tag 2.11
$$
Let us estimate the second moment of the fourth and fifth summands
in (2.11) (the first moment of the second and third summands can
be estimated in the same way). Using the fact that the processes
$\{\eta_t\}$ and $\{\widetilde{w}(t)\}$ are independent we obtain
$$
{\Bbb E}\left|\int_{\eta_t}^{\eta_t^n}\widetilde{\beta}(x_n^S
(\rho_u^n))d\widetilde{w}(u)\right|^2={\Bbb E}\left[\left.{\Bbb E}
\left|\int_{\eta_t}^{\eta_t^n}\widetilde{\beta}(x_n^S(\rho_u^n))
d\widetilde{w}(u)\right|^2\right/{\Cal F}^w_\infty\right]=
$$
$$
={\Bbb E}\left[\left.{\Bbb E}\int_{\eta_t}^{\eta_t^n}\left\|
\widetilde{\beta}(x_n^S(\rho_u^n))\right\|^2du\right/{\Cal
F}^w_\infty \right]\leq K^2{\Bbb E}|\eta_t-\eta_t^n|\to 0,\quad
n\to+\infty.
$$
Calculating in the same way the second moment we get
$$
{\Bbb E}\left|\int_0^{\eta_t}\left(\widetilde{\beta}
(x_n^S(\rho_u^n))-\widetilde{\beta}(x^S(\rho_u))\right)d\widetilde{w}
(u)\right|^2={\Bbb E}\int_0^{\eta_t}\left\|\widetilde{\beta}
(x_n^S(\rho_u^n))-\widetilde{\beta}(x^S(\rho_u))\right\|^2 du\to 0
$$
as $n\to+\infty.$ We take into account the following arguments:
$\widetilde{\beta}(x_n^S(\rho_u^n))\mathop{\longrightarrow}\limits^P
\widetilde{\beta}(x^S(\rho_u)),$ $n\to+\infty$ (from Lemma 4) and
$\left\|\widetilde{\beta}(x_n^S(\rho_u^n))-\widetilde{\beta}
(x^S(\rho_u))\right\|\leq 2K.$

For estimating the last summand in (2.11) and the first summand in
(2.10) we need the following result.

\proclaim{Lemma 5} For all $t\geq 0: \rho^n_{\eta_t^n}
\mathop{\longrightarrow}\limits^P\rho_{\eta_t}$ when $n\to +\infty.$
\endproclaim

\demo{Proof} Similarly to Lemma 2 we get for $\varepsilon>0:$
$$
{\Bbb P}\left\{\rho_{\eta_t}-\rho^n_{\eta_t^n}>\varepsilon\right\}
=\sum_{k\geq 0}{\Bbb P}\left\{ \rho_{\eta_t}>\varepsilon+
\rho^n_{\eta_t^n},\rho_{\eta_t}\in\left({(k+1)\varepsilon\over2},{(k+
2)\varepsilon\over2}\right]\right\}\leq
$$
$$
\leq\sum_{k\geq 0}{\Bbb P}\left\{\rho^n_{\eta_t^n}<
{k\varepsilon\over2},\rho_{\eta_t}>{(k+1)\varepsilon\over 2}\right\}=
\sum_{k\geq 0}{\Bbb P}\left\{\eta_t^n\leq\eta^n_{{k\varepsilon\over
2}},\eta_t>\eta_{{(k+1)\varepsilon\over 2}}\right\}.\tag 2.12
$$
Consider the $k-$th summand in (2.12). For some $\delta>0$ we have
$$
{\Bbb P}\left\{\eta_t^n\leq\eta^n_{{k\varepsilon\over2}},\eta_t>
\eta_{{(k+1)\varepsilon\over 2}}\right\}\leq {\Bbb P}\left\{
\eta_t^n\leq\eta^n_{{k\varepsilon\over2}},\eta_t-\delta\geq
\eta_{{(k+1)\varepsilon\over2}}\right\}+
$$
$$
+{\Bbb P}\left\{\eta_t-\delta<\eta_{{(k+1)\varepsilon\over2}}<\eta_t
\right\}\leq {\Bbb P}\left\{\eta_t^n\leq\eta^n_{{k\varepsilon
\over2}},\eta_t-\delta\geq\eta_{{(k+1)\varepsilon\over2}},|\eta_t^n-
\eta_t|\leq{\delta\over2}\right\}+
$$
$$
+{\Bbb P}\left\{|\eta_t^n-\eta_t|>{\delta\over2}\right\}+ {\Bbb P}
\left\{\eta_t-\delta<\eta_{{(k+1)\varepsilon\over2}}<\eta_t\right\}
\leq {\Bbb P}\left\{|\eta_{{(k+1)\varepsilon\over2}}-\eta_{{(k+1)
\varepsilon\over2}}^n|>{\delta\over2}\right\}+
$$
$$
+{\Bbb P}\left\{|\eta_t^n-\eta_t|>{\delta\over2}\right\}+{\Bbb P}
\left\{\eta_t-\delta<\eta_{{(k+1)\varepsilon\over2}}<\eta_t\right\}.
\tag 2.13
$$
In the last inequality we use that
$$
\left\{\eta_t^n\leq\eta^n_{{k\varepsilon\over2}},\eta_t-\delta\geq
\eta_{{(k+1)\varepsilon\over2}},|\eta_t^n-\eta_t|\leq{\delta\over2}
\right\}\subseteq\left\{|\eta_{{(k+1)\varepsilon\over2}}-\eta_{{(k+1)
\varepsilon\over2}}^n|>{\delta\over2}\right\}.
$$

Proceeding $n\to+\infty$ and using Corollary 2 (when $q\in(-1,0)\cup
(0,1)$) or Remark 1 (when $|q|=1$) or Remark 2 (when $q=0$) we see
that the first and second summands are equal to $0.$ Consider the
third summand. Note that $\{\eta_t\}$ is additive functional of the
Markov process $\{x^\nu(t)\}$ and for all $a>0:$ ${\Bbb P}_x\{
\eta_t<a\}\leq {a\sqrt{2}\over\sqrt{\pi t}}$ (${\Bbb P}_x$ is
standard notation for ${\Bbb P}\{\cdot/x_\nu(0)=x\}$). Also ${\Bbb
P}\left\{\eta_t-\delta<\eta_{{(k+1)\varepsilon\over2}}<\eta_t\right\}
=0$ when $t\geq{(k+1)\varepsilon\over2}.$ Therefore for all $t<{(k+1)
\varepsilon\over2}$ we get
$$
{\Bbb P}\left\{\eta_t-\delta<\eta_{{(k+1)\varepsilon\over2}}<\eta_t
\right\}={\Bbb E}{\Bbb P}_{x^\nu\left({(k+1)\varepsilon\over2}
\right)}\left\{0<\eta_{t-{(k+1)\varepsilon\over2}}<\delta\right\}\leq
{\delta\sqrt{2}\over\sqrt{\pi (t-{(k+1)\varepsilon\over2})}}\to 0
$$
when $\delta\to 0.$ In the same way one can prove that ${\Bbb P}
\left\{\rho^n_{\eta_t^n}-\rho_{\eta_t}>\varepsilon\right\}\to 0,
n\to+\infty.$ The lemma is proved.

Due to Lemma 5 for the last summand in (2.11) and the first summand
in (2.10), for some $N>0$ we have
$$
{\Bbb P}\{|w(\rho^n_{\eta_t^n})-w(\rho_{\eta_t})|>\varepsilon\}\leq
{\Bbb P}\{\|w\|_{\,\hbox{{\rm Hol}}\,_{1/4}([0,t])}|\rho^n_{\eta_t^n}
-\rho_{\eta_t}|^{1/4}>\varepsilon,\|w\|_{\,\hbox{{\rm Hol}} \,_{1/4}
([0,t])}<N\}+
$$
$$
+{\Bbb P}\{\|w\|_{\,\hbox{{\rm Hol}}\,_{1/4}([0,t])}\geq N\}\leq{\Bbb
P}\left\{|\rho^n_{\eta_t^n}-\rho_{\eta_t}|^{1/4}>{\varepsilon\over
N}\right\}+{\Bbb P}\{\|w\|_{\,\hbox{{\rm Hol}}\,_{1/4}([0,t] )}\geq
N\}.\tag 2.14
$$
Let $n\to+\infty,$ then the first summand in (2.14) tends to $0.$
Then let $N\to+\infty.$ We obtain that the last summand in (2.12) and
the first summand in (2.10) tend to 0. This completes the proof of
Theorem 2.\enddemo

\proclaim{Corollary 3} It follows from Theorem 2 that solution of
(2.1), considered as a random function on $\Re^+\times\Re^d,$ has
a measurable modification.\endproclaim

\end